\documentclass[11]{article}
\usepackage[russian]{babel}
\usepackage{amssymb}
\usepackage{inputenc}
\usepackage{graphicx}
\begin{document}
\begin{center}
\textbf{\Large CALCULATION OF THE NORM OF THE ERROR FUNCTIONAL OF
OPTIMAL QUADRATURE FORMULAS\\ IN THE SPACE $W_2^{(2,1)}(0,1)$ }\\
\textbf{A.R.Hayotov}
\end{center}

\textbf{Abstract.} In this paper in the space $W_2^{(2,1)}(0,1)$
square of the norm of the error functional of a optimal
quadrature formula is calculated.

\large

Many well known mathematicians be occupied with construction of
optimal formulas for approximate integration. Full bibliography by
this direction can be found in books [1,2,3,4].

In works [5,6] in space $W_2^{(m,m-1)}(0,1)$ was considered
problem of construction of the optimal quadrature formulas of the
form
$$
\int\limits_0^1p(x)\varphi(x)dx\cong \sum\limits_{\beta=0}^N
C_{\beta}\varphi(x_\beta)\eqno (1)
$$
with error functional
$$
\ell(x)=p(x)\varepsilon_{[0,1]}(x)- \sum\limits_{\beta=0}^N
C_{\beta}\delta(x-x_\beta),\eqno (2)
$$
where $C_{\beta}$ and $x_{\beta}\in [0,1]$
($\beta=\overline{0,N}$) are called \emph{coefficients} and
\emph{nodes} of the quadrature formula (1) respectively, $p(x)$ is
a weight function, $\varepsilon_{[0,1]}(x)$ is the characteristic
function of the interval [0,1], $\delta(x)$ is Dirac delta
function, and $\varphi(x)$ is such a function, that is contained
in Hilbert space $W_2^{(m,m-1)}(0,1)$ norm of functions in which
is defined by formula
$$
\|\varphi(x)|W_2^{(m,m-1)}(0,1)\|=\left\{\int\limits_0^1\left(\varphi^{(m)}(x)+
\varphi^{(m-1)}(x)\right)^2dx\right\}^{1/2}.
$$

Note that in work [5] was found the extremal function and with its
help following representation of square of the norm of the error
functional (2) was obtained:
$$
\|\ell(x)\|^2_{W_2^{(m,m-1)*}(0,1)}=(-1)^m\Bigg[\sum\limits_{\beta=0}^N\sum\limits_{\beta'=0}^N
C_{\beta}C_{\beta'}\psi_m(x_{\beta}-x_{\beta'})-
$$
$$
-2\sum\limits_{\beta=0}^NC_{\beta}\int\limits_0^1p(x)\psi_m(x-x_{\beta})dx+
\int\limits_0^1\int\limits_0^1p(x)p(y)\psi_m(x-y)dxdy\Bigg],\eqno
(3)
$$
where
$$
\psi_m(x)=\frac{\mathrm{sign}x}{2}\left(\frac{e^x-e^{-x}}{2}-
\sum\limits_{k=1}^{m-1}\frac{x^{2k-1}}{(2k-1)!}\right). \eqno (4)
$$
Moreover, the error functional (2), as shown in [5], satisfies
following orthogonality conditions
$$
(\ell(x),x^{\alpha})=0,\ \ \alpha=\overline{0,m-2},\eqno (5)
$$
$$
(\ell(x),e^{-x})=0. \eqno (6)
$$

The norm (3) of the error functional $\ell(x)$ is many-dimensional
function of the coefficients $C_{\beta}$ ($\beta=\overline{0,N}$).
Since error of the quadrature formula (1) is estimated from above
by the norm of the error functional  $\ell(x)$ in conjugate space,
then in order to construct the optimal quadrature formula of the
form (1) it is required to minimize, taking account of conditions
(5) and (6), square of the norm (3) by coefficients $C_{\beta}$
when the nodes $x_{\beta}$ are fixed, i.e. we need find condition
minimum of square of the error functional norm in conditions (5),
(6).

Further in [5], applying the method of Lagrange undetermined
factors, for finding of minimum of square of the error functional
norm following discrete system of Wiener-Hopf type was obtained
$$
\sum_{\gamma=0}^NC_{\gamma}\psi_m(x_{\beta}-x_{\gamma})+P_{m-2}(x_{\beta})+de^{-x_{\beta}}=\int\limits_0^1
p(x)\psi_m(x-x_{\beta})dx,\ \beta=\overline{0,N},\eqno(7)
$$
$$
\sum_{\gamma=0}^N
C_{\gamma}x_{\gamma}^{\alpha}=\int\limits_0^1p(x)x^{\alpha}dx,\
\alpha=\overline{0,m-2}, \eqno (8)
$$
$$
\sum_{\gamma=0}^N
C_{\gamma}e^{-x_{\gamma}}=\int\limits_0^1p(x)e^{-x}dx, \eqno (9)
$$
where $\psi_m(x)$ is defined by formula (4), $P_{m-2}(x)$ is
unknown polynomial of degree $m-2$, $d$ is unknown constant. And
also existence and uniqueness of the solution of this system were
proved.

In [6] when $x_{\beta}=h\beta$, $\beta=\overline{0,N},$
$h=\frac{1}{N}$, $N=1,2,...$ the system (7)-(9) was solved and it
was found the analytical representations of the coefficients
$C_{\beta}$ ($\beta=\overline{0,N}$). The coefficients for which
minimum of square of the error functional norm is attained are
called by \emph{optimal}. In particular, for $m=2$, $p(x)=1$
following theorem was proved

\textbf{Theorem 1.} {\it The coefficients of optimal quadrature
formulas of the form (1) in the space $W_2^{(2,1)}(0,1)$ when
$p(x)=1$ have following view: }
$$
C_{\beta}=\left\{
\begin{array}{ll}
\frac{e^h-1-h}{e^h-1}-K(h)(\lambda_1-\lambda_1^N), &\beta=0 \\
h-K(h)\left[(\lambda_1-e^h)\lambda_1^{\beta}+(\lambda_1
e^h-1)\lambda_1^{N-\beta}\right],& \beta=\overline{1,N-1},\\
\frac{he^h-e^h+1}{e^h-1}-K(h)(\lambda_1-\lambda_1^N)e^h, &\beta=N, \\
\end{array}
\right.\eqno (10)
$$
\emph{where}
$$
K(h)=\frac{(2e^h-2-he^h-h)(\lambda_1-1)}{2(e^h-1)^2(\lambda_1+\lambda_1^{N+1})},
\eqno (11)
$$
$$
\lambda_1=\frac{h(e^{2h}+1)-e^{2h}+1-(e^h-1)\sqrt{h^2(e^h+1)^2+2h(1-e^h)}}
{1-e^{2h}+2he^h}. \eqno (12)
$$

Present paper is direct continuation of works [5,6]. Aim of given
work is with using theorem 1 to calculate the square of the norm
of the error functional $\ell(x)$ in the space
$W_2^{(2,1)*}(0,1)$.

Following is take placed

\textbf{Теорема 2.} {\it For square of the norm (3) of the
functional (2) of optimal quadrature formula (1) when $p(x)=1$,
$x_{\beta}=h\beta$ in the space $W_2^{(2,1)*}(0,1)$ following
equality is valid
$$
\|\ell(x)\|^2_{W_2^{(2,1)*}(0,1)}=\frac{h^2}{12}+
\frac{h(2-e^h-3e^{2h})+4+2e^h+6e^{2h}}{4(1-e^h)^2}+
$$
$$
+K(h)\bigg[\frac{(\lambda_1^N+\lambda_1^2)(1+e^h)-(\lambda_1^{N+1}+\lambda_1)
(1+2e^h)}{2(1-\lambda_1)}+
$$
$$
+\frac{h^2(\lambda_1^2+\lambda_1)(\lambda_1^N-1)(1+e^h)}
{2(1-\lambda_1)^2}+
$$
$$
+\frac{(\lambda_1-e^h)^2(\lambda_1^N-\lambda_1e^h)-
(1-\lambda_1e^h)^2(\lambda_1-\lambda_1^Ne^h)}{2(1-\lambda_1e^h)(\lambda_1-e^h)}\bigg],
$$
where $K(h)$ and $\lambda_1$ are determined by (11) and (12)
respectively. }

\textbf{Proof.} The system (7)-(9) for $m=2,$ $p(x)=1$,
$x_{\beta}=h\beta$ have form
$$
\sum\limits_{\gamma=0}^NC_{\gamma}\psi_2(h\beta-h\gamma)+P_0(h\beta)+d\
e^{-h\beta}=\int\limits_0^1\psi_2(x-h\beta)dx,\eqno (13)
$$
$$
\sum\limits_{\gamma=0}^NC_{\gamma}=1,\eqno (14)
$$
$$
\sum\limits_{\gamma=0}^NC_{\gamma}e^{-h\gamma}=1-e^{-1}.\eqno (15)
$$
Then for square of the norm (3) of the functional $\ell(x)$ we
obtain
$$
\|\ell(x)\|^2=\sum\limits_{\beta=0}^NC_{\beta}\left(\sum\limits_{\beta'=0}^NC_{\beta'}
\psi_2(h\beta-h\beta')-\int\limits_0^1\psi_2(x-h\beta)dx\right)-
$$
$$
-\sum\limits_{\beta=0}^NC_{\beta}\int\limits_0^1\psi_2(x-h\beta)dx+
\int\limits_0^1\int\limits_0^1\psi_2(x-y)dxdy.
$$
Hence, taking into account (13), we have
$$
\|\ell(x)\|^2=-\sum\limits_{\beta=0}^NC_{\beta}\left(P_0(h\beta)+d\
e^{-h\beta}\right)-
$$
$$
-\sum\limits_{\beta=0}^NC_{\beta}\int\limits_0^1\psi_2(x-h\beta)dx+
\int\limits_0^1\int\limits_0^1\psi_2(x-y)dxdy.\eqno (16)
$$
From here, taking account of (4), for integrals of (16) we get
$$
\int\limits_0^1\psi_2(x-h\beta)dx=\frac{e^{h\beta}+e^{-h\beta}+e^{1-h\beta}+e^{h\beta-1}-4}
{4}-\frac{(h\beta)^2+(1-h\beta)^2}{4},\eqno (17)
$$
$$
\int\limits_0^1\int\limits_0^1\psi_2(x-y)dxdy=\frac{e^2-1}{2e}-\frac{7}{6}.
\eqno (18)
$$

In representation (16) of the error functional norm polynomial
$P_0(h\beta)=b_0$ and constant $d$ are unknowns. For
$\|\ell(x)\|^2$, using (14), (15), (17), (18), we have
$$
\|\ell(x)\|^2=-b_0+\frac{1-e}{e}d-\frac{e+1}{4e}\sum\limits_{\beta=0}^NC_{\beta}e^{h\beta}
-\frac{1+e}{4}(1-e^{-1})+\frac{5}{4}+
$$
$$
+\frac{1}{2}\sum\limits_{\beta=0}^NC_{\beta}(h\beta)^2-\frac{1}{2}\sum\limits_{\beta=0}^NC_{\beta}
(h\beta)+\frac{e^2-1}{2e}-\frac{7}{6}.\eqno (19)
$$

The equality (13) take placed in any $h\beta$ when
$\beta=\overline{0,N}$, $h=\frac1N$, i.e. is identity by powers
$h\beta$ and $e^{h\beta}$, $e^{-h\beta}$. Equating coefficients of
the left and the right sides of (13) in front of $e^{-h\beta}$ and
constant term, by using theorem 1, equalities (4), (14), (15),
(17), for $d$ and $b_0$ we get
$$
d=\frac{C_0}{2}+\frac12\left(\frac{he^h}{1-e^h}+a_1\frac{\lambda_1e^h}{1-\lambda_1e^h}
+b_1\frac{\lambda_1^Ne^h}{\lambda_1-e^h}\right)-\frac{1}{4}\sum\limits_{\gamma=0}^NC_\gamma
e^{h\gamma}+\frac{1+e}{4},
$$
$$
-b_0=\frac{h(1+e^h)}{2(1-e^h)}+ha_1\frac{\lambda_1}
{(1-\lambda_1)^2}+hb_1\frac{\lambda_1^{N+1}}{(1-\lambda_1)^2}-
\frac12\sum\limits_{\gamma=0}^NC_\gamma(h\gamma)+\frac54,
$$
where $a_1=K(h)(e^h-\lambda_1)$, $b_1=K(h)(1-\lambda_1e^h)$.\\
Then from (19), taking account of $d$ and $b_0$, using identities
$$
\sum\limits_{\gamma=1}^{N-1}\lambda^\gamma\gamma=\frac{\lambda-\lambda^{N+1}-
N\lambda^N(1-\lambda)}{(1-\lambda)^2},
$$
$$
\sum\limits_{\gamma=1}^{N-1}\lambda^\gamma\gamma^2=\frac{\lambda^N(\lambda^2+\lambda+
N^2(1-\lambda)^2+2N(\lambda-\lambda^2))}{(\lambda-1)^3}-\frac{\lambda^2+\lambda}{(\lambda-1)^3}
$$
and theorem 1, after simplification we get the statement of the
theorem 2.
\newpage

\textbf{References}
\begin{enumerate}
\item Krylov V.I. Approximate Calculation of Integrals. - M.: Nauka, 1967. - 500p.
\item Nikolskii S.M. Quadrature Formulas. - M.: Nauka, 1988. - 256 p.
\item Sobolev S.L. Introduction to the Theory of Cubature Formulas. -M.: Nauka, 1974.
- 808p.
\item
Sobolev S.L., Vaskevich V.L. Cubature Formulas. -Novosibirsk,
Institute of Mathematics SB of RAS, 1996, -484p.
\item Shadimetov Kh.M., Hayotov A.R. Weight Optimal Quadrature Formulas
in the Space\\ $W_2^{(m,m-1)}(0,1)$. Uzbek Mathematical Journal.
2002. \No 3-4, pp.92-103.
\item Shadimetov Kh.M., Hayotov A.R. Calculation of the Coefficients
of Optimal Quadrature Formulas in the $W_2^{(m,m-1)}(0,1)$. Uzbek
Mathematical Journal. 2004. \No 3. pp.67-82.\\
(arXiv.0810.5421v1 [math.NA]).
\end{enumerate}

\noindent
Abdullo Rakhmonovich Hayotov\\
Institute of Mathematics and Information Technologies\\
Uzbek Academy of Sciences\\
F.Hodjaev str, 29\\
Tashkent, 100125\\
Uzbekistan\\
\textit{E-mail:} abdullo\_hayotov@mail.ru, hayotov@mail.ru.

\end{document}